\documentclass[authoryear]{elsarticle}
\usepackage{amsfonts, amsmath, amssymb, amsbsy, amsthm, graphicx, calrsfs, color}

\numberwithin{equation}{section}

\theoremstyle{plain}
\newtheorem{theorem}{Theorem}[section]
\newtheorem{proposition}[theorem]{Proposition}

\theoremstyle{remark}
\newtheorem{remark}[theorem]{Remark}

\newcommand{\Acfg}{\hat{A}_n^{\mathrm{CFG}}}
\newcommand{\Azwp}{\hat{A}_n^{\mathrm{ZWP}}}
\newcommand{\Adh}{\hat{A}_n^{\mathrm{D}}}
\newcommand{\Aht}{\hat{A}_n^{\mathrm{HT}}}
\newcommand{\Anaive}{\hat{A}_n}
\newcommand{\Apick}{\hat{A}_n^{\mathrm{P}}}
\newcommand{\Acfgad}{\hat{A}_{n,\mathrm{ad}}^{\mathrm{CFG}}}
\newcommand{\Aols}{\hat{A}_n^{\mathrm{OLS}}}
\newcommand{\sigmaols}{\hat{\sigma}_{n,\mathrm{OLS}}^2}

\newcommand{\opt}{^{\mathrm{opt}}}
\newcommand{\etaopt}{\eta_{\mathrm{opt}}}
\newcommand{\Acfgopt}{\hat{A}_{n,\mathrm{opt}}^{\mathrm{CFG}}}

\newcommand{\mv}[1]{\boldsymbol{#1}}
\newcommand{\1}{\boldsymbol{1}}
\newcommand{\simplex}{{\Delta_p}}
\newcommand{\RR}{\mathbb{R}}
\newcommand{\cont}{{\cal C}}
\newcommand{\bounded}{\ell^\infty}
\newcommand{\Fclass}{{\cal F}}
\newcommand{\dto}{\rightsquigarrow}

\renewcommand{\le}{\leqslant}
\renewcommand{\ge}{\geqslant}
\newcommand{\cov}{\operatorname{cov}}
\newcommand{\var}{\operatorname{var}}

\begin{document}
\begin{frontmatter}
\title{Nonparametric estimation of an extreme-value copula in arbitrary dimensions}
\author{Gordon Gudendorf\fnref{ARC-PAI}}
\ead{gordon.gudendorf@uclouvain.be}
\author{Johan Segers\corref{cor}\fnref{ARC-PAI}}
\ead{johan.segers@uclouvain.be}
\address{Institut de Statistique, Universit\'{e} catholique de Louvain, Voie du Roman Pays 20, B-1348 Louvain-la-Neuve, Belgium}
\cortext[cor]{Corresponding author}
\fntext[ARC-PAI]{Research supported by IAP research network grant nr.\ P6/03 of the Belgian government (Belgian Science Policy) and by contract nr.\ 07/12/002 of the Projet d'Actions de Recherche Concert\'ees of the Communaut\'e fran\c{c}aise de Belgique, granted by the Acad\'emie universitaire Louvain.}

\begin{abstract}
Inference on an extreme-value copula usually proceeds via its Pickands dependence function, which is a convex function on the unit simplex satisfying certain inequality constraints. In the setting of an iid random sample from a multivariate distribution with known margins and unknown extreme-value copula, an extension of the Cap\'era\`a--Foug\`eres--Genest estimator was introduced by D.~Zhang, M.~T.~Wells and L.~Peng [Journal of Multivariate Analysis 99 (2008) 577--588]. The joint asymptotic distribution of the estimator as a random function on the simplex was not provided. Moreover, implementation of the estimator requires the choice of a number of weight functions on the simplex, the issue of their optimal selection being left unresolved. 

A new, simplified representation of the CFG-estimator combined with standard empirical process theory provides the means to uncover its asymptotic distribution in the space of continuous, real-valued functions on the simplex. Moreover, the ordinary least-squares estimator of the intercept in a certain linear regression model provides an adaptive version of the CFG-estimator whose asymptotic behavior is the same as if the variance-minimizing weight functions were used. As illustrated in a simulation study, the gain in efficiency can be quite sizeable.
\end{abstract}



\end{frontmatter}

\section{Introduction}
\label{S:intro}

Let $\mv{X}_i = (X_{i1}, \ldots, X_{ip})$, $i \in \{1, \ldots, n\}$, be iid random vectors from a $p$-variate, continuous distribution function $F$ with multivariate extreme-value copula $C$: for $\mv{u} \in (0, 1]^p \setminus \{(1, \ldots, 1)\}$, denoting the margins of $F$ by $F_1, \ldots, F_p$,
\begin{multline}
\label{E:MEVC}
  C(\mv{u})
  = P \bigl( F_1(X_{i1}) \le u_1, \ldots, F_p(X_{ip}) \le u_p \bigr)
  = \exp \{ - |\mv{y}| \, A( \mv{y} / |\mv{y}| ) \} \\
  \qquad \text{where $y_j = - \log u_j$ and $|\mv{y}| = |y_1| + \cdots + |y_p|$}.
\end{multline}
The function $A$, whose domain is $\simplex = \{ \mv{w} \in [0, 1]^p : w_1 + \cdots + w_p = 1 \}$, is called the Pickands dependence function of $C$, after \cite{Pickands81}. 

Multivariate extreme-value copulas arise as the limits of copulas of vectors of component-wise maxima of independent random samples \citep{Deheuvels84, Galambos87}. As a consequence, they coincide with the class of copulas of multivariate extreme-value or max-stable distributions. Therefore, they provide models for dependence between extreme values that allow extrapolation beyond the support of the sample. It is then of interest to estimate the Pickands dependence function $A$.

A necessary condition for $C$ in \eqref{E:MEVC} to be a copula is that $A$ is convex and satisfies $\max(w_1, \ldots, w_p) \le A(\mv{w}) \le 1$ for all $\mv{w} \in \simplex$; in the bivariate case, this is also sufficient. In general, $A$ should admit an integral representation in terms of a spectral measure. Some other properties of Pickands dependence functions are studied in \cite{Obretenov91} and \cite{FRE05}. The upshot of all this is that the class of Pickands dependence functions is infinite-dimensional. This warrants the use of nonparametric methods. 

Whereas most papers hitherto concentrated on the bivariate case, a nonparametric estimator for general multivariate Pickands dependence functions was introduced in \cite{ZWP08}. This estimator is in fact a multivariate generalization of the one by Cap\'era\`a--Foug\`eres--Genest \citep{CFG97}. The estimator was shown to be uniformly consistent and pointwise asymptotically normal. However, the joint asymptotic distribution of the estimator as a random function on $\simplex$ was not provided. Moreover, implementation of the estimator requires the choice of $p$ weight functions $\lambda_j$ on $\simplex$, the issue of their optimal selection being left unresolved.

Using a simplified representation of the above-mentioned estimator, we are able to uncover its asymptotic distribution in the space $\cont(\simplex)$ of continuous, real-valued functions on $\Delta_p$. Moreover, we give explicit expressions for the weight functions $\lambda_j$ that minimize the pointwise asymptotic variance of the estimator. These optimal weight functions depend on the unknown distribution. We show that the CFG-estimator with estimated variance-minimizing weight functions can be implemented as the intercept estimator in a certain linear regression model via ordinary least squares. The OLS-estimator is data-adaptive in the sense that the asymptotic distribution is the same as if the optimal weight functions were used. In a simulation study, the gain in efficiency is shown to be quite sizeable.

As in \cite{ZWP08}, the setting here is that of a random sample from a distribution whose margins are known and whose copula is an extreme-value copula. It would be worthwhile to extend this to the case of unknown margins \citep{GP08, GS09} and the case that the copula of $F$ is merely in the domain of attraction of an extreme-value copula \citep{CF00, ES09}.

The outline of our paper is as follows. The CFG-estimator is introduced in the next section, including its simplified representation and asymptotic distribution. The variance-minimizing weight functions are computed in Section~\ref{S:OLS} together with an adaptive estimator based on ordinary least squares in a linear regression framework. Section~\ref{S:simul} reports on a simulation study. The proofs of the results in Sections~\ref{S:CFG} and~\ref{S:OLS} are deferred to Appendices~A and~B, respectively.

\section{CFG-estimator and variants}
\label{S:CFG}

Let $\mv{X}_i = (X_{i1}, \ldots, X_{ip})$, $i \in \{1, \ldots, n\}$, be iid random vectors from a $p$-variate, continuous distribution function $F$ with multivariate extreme-value copula $C$ and Pickands dependence function $A$ as in \eqref{E:MEVC}. Let $F_1, \ldots, F_p$ be the marginal distribution functions of $F$. Put $\mv{Y}_i = (Y_{i1}, \ldots, Y_{ip})$ where
\begin{equation}
\label{E:Yi}
  Y_{ij} = - \log F_j(X_{ij})
\end{equation}
for $i \in \{1, \ldots, n\}$ and $j \in \{1, \ldots, p\}$. The marginal distributions of the random variables $Y_{ij}$ are standard exponential. The random vectors $\mv{Y}_1, \ldots, \mv{Y}_p$ are iid with common joint survivor function
\[
  P(Y_{i1} > y_1, \ldots, Y_{ip} > y_p)
  = C(e^{-y_1}, \ldots, e^{-y_p}) 
  = \exp \{ - |\mv{y}| \, A(\mv{y} / |\mv{y}|) \},
\]
for $\mv{y} \in [0, \infty)^p \setminus \{(0, \ldots, 0)\}$, where $|\mv{y}| = |y_1| + \cdots + |y_p|$. Put
\begin{equation}
\label{E:xi}
  \xi_i(\mv{w}) = \bigwedge_{j=1}^p \frac{Y_{ij}}{w_j}, \qquad \mv{w} \in \simplex, \; i \in \{1, \ldots, n\},
\end{equation}
with `$\wedge$' denoting minimum and with the obvious convention for division by zero; in particular, $\xi_i(\mv{e}_j) = Y_{ij}$ for the $p$ standard unit vectors $\mv{e}_1, \ldots, \mv{e}_p$ in $\RR^p$. For $\mv{w} \in \simplex$ and $x > 0$, we have
\begin{equation}
\label{E:xiexp}
  P \bigl( \xi_i(\mv{w}) > x \bigr) = P(Y_{i1} > w_1 x, \ldots, Y_{ip} > w_p x) = \exp \{ - x \, A(\mv{w}) \}.
\end{equation}
Hence the random variables $\xi_1(\mv{w}), \ldots, \xi_n(\mv{w})$ constitute an independent random sample from the exponential distribution with mean $1/A(\mv{w})$. It follows that the distribution of $- \log \xi_i(\mv{w})$ is Gumbel with location parameter $\log A(\mv{w})$, whence
\begin{equation}
\label{E:Elogxi}
  E [ - \log \xi_i(\mv{w}) ] = \log A(\mv{w}) + \gamma,
\end{equation}
the Euler--Mascheroni constant $\gamma = - \Gamma'(1) = 0.5772\ldots$ being the mean of the standard Gumbel distribution. This suggests the naive estimator
\begin{equation}
\label{E:Anaive}
  \log \Anaive(\mv{w}) = - \frac{1}{n} \sum_{i=1}^n \log \xi_i(\mv{w}) - \gamma, \qquad \mv{w} \in \simplex.
\end{equation}
The naive estimator is itself not a valid Pickands dependence function. For instance, it does not verify the vertex constraints $A(\mv{e}_j) = 1$ for all $j \in \{1, \ldots, p\}$. A simple way to at least remedy this defect is by putting
\begin{equation}
\label{E:Acfg}
  \log \Acfg(\mv{w}) = \log \hat{A}_n(\mv{w}) - \sum_{j=1}^p \lambda_j(\mv{w}) \, \log \hat{A}_n(\mv{e}_j), \qquad \mv{w} \in \simplex,
\end{equation}
where $\lambda_1, \ldots, \lambda_p : \simplex \to \RR$ are continuous functions verifying $\lambda_j(\mv{e}_k) = \delta_{jk}$ for all $j,k \in \{1, \ldots, p\}$. Continuity of the functions $\lambda_j$ is assumed merely to ensure that the resulting estimator is a continuous function of $\mv{w}$ as well. 

The superscript `CFG' refers to the bivarate estimator by Cap\'era\`a--Foug\`eres--Genest in \cite{CFG97}, generalized to the multivariate case in \cite{ZWP08}. Actually, the original definition in \cite{ZWP08} is
\begin{equation}
\label{E:Azwp}
  \log \Azwp(\mv{w}) = \sum_{j=1}^p \lambda_j(\mv{w}) \; \int_0^{1-w_j} \frac{n^{-1} \sum_{i=1}^n \1\{Z_{ij}(\mv{w}) \le z\} - z}{z(1-z)} \, dz,
\end{equation}
where, with $Y_{ij}$ as in \eqref{E:Yi},
\[
  Z_{ij}(\mv{w}) = \frac{\bigwedge_{k : k \neq j} \frac{Y_{ik}}{w_k}}{\frac{Y_{ij}}{1-w_j} + \bigwedge_{k : k \neq j} \frac{Y_{ik}}{w_k}},
  \qquad \mv{w} \in \simplex.
\]
Moreover, in \eqref{E:Azwp}, the weight functions $\lambda_j$ are supposed to be nonnegative and to satisfy the additional constraint
\begin{equation}
\label{E:ZWP}
  \sum_{j=1}^p \lambda_j(\mv{w}) = 1, \qquad \mv{w} \in \simplex.
\end{equation}
However, if \eqref{E:ZWP} holds, then actually the two estimators coincide, that is,
\begin{equation}
\label{E:ZWPisCFG}
  \Azwp(\mv{w}) = \Acfg(\mv{w}), \qquad \mv{w} \in \simplex.
\end{equation}
The proof of \eqref{E:ZWPisCFG} is essentially the same as the one in \cite{S07Ahs} for the bivariate case, the key being that the integrals in~\eqref{E:Azwp} can be solved:
\begin{multline*}
  \int_0^{1-w_j} \frac{\1\{Z_{ij}(\mv{w}) \le z\} - z}{z(1-z)} \, dz \\
  = \log [ 1 - \{(1-w_j) \wedge Z_{ij}(\mv{w})\} ] + \log(1 - w_j) - \log \{ (1-w_j) \wedge Z_{ij}(\mv{w}) \} \\
  = \log Y_{ij} - \log \xi_i(\mv{w}).
\end{multline*}
In our representation \eqref{E:Acfg}, however, there is no reason whatsoever to restrict the weight functions to satisfy~\eqref{E:ZWP}.

The asymptotics of the naive estimator and the CFG-estimator follow from standard empirical process theory as presented for instance in \cite{VW96} and \cite{VdV98}. Let $\cont(\simplex)$ denote the Banach space of continuous functions from $\simplex$ into $\RR$ equipped with the supremum norm. Convergence in distribution is denoted by the arrow `$\dto$'.

\begin{proposition}[Naive estimator]
\label{P:naive}
Let $\mv{X}_i = (X_{i1}, \ldots, X_{ip})$, $i \in \{1, \ldots, n\}$, be iid random variables from a $p$-variate, continuous distribution function $F$ with multivariate extreme-value copula $C$ and Pickands dependence function $A$. The naive estimator $\Anaive$ in \eqref{E:Anaive} satisfies
\begin{equation}
\label{E:naive:cons}
  \sup_{\mv{w} \in \simplex} |\Anaive(\mv{w}) - A(\mv{w})| \to 0, \qquad n \to \infty, \qquad \text{almost surely},
\end{equation}
and in $\cont(\simplex)$,
\begin{equation}
\label{E:naive:CLT}
  \sqrt{n} (\Anaive - A) \dto A \, \zeta, \qquad n \to \infty,
\end{equation}
where $\zeta$ is a centered Gaussian process with covariance function
\begin{equation}
\label{E:cov:zeta}
  \cov \bigl( \zeta(\mv{v}), \zeta(\mv{w}) \bigr) = \cov \bigl( - \log \xi_i(\mv{v}), \, - \log \xi_i(\mv{w}) \bigr),
  \qquad \mv{v}, \mv{w} \in \simplex,
\end{equation}
with $\xi_i(\,\cdot\,)$ as in \eqref{E:xi}.
\end{proposition}

\begin{theorem}[CFG-estimator]
\label{T:CFG}
If, in addition to the assumptions in Proposition~\ref{P:naive}, the functions $\lambda_1, \ldots, \lambda_p : \simplex \to \RR$ are continuous, then
\begin{equation}
\label{E:CFG:cons}  
\sup_{\mv{w} \in \simplex} |\Acfg(\mv{w}) - A(\mv{w})| \to 0, \qquad n \to \infty, \qquad \text{almost surely},
\end{equation}
and in $\cont(\simplex)$,
\begin{equation}
\label{E:CFG:CLT}
  \sqrt{n} (\Acfg - A) \dto A \, \eta, \qquad n \to \infty,
\end{equation}
where $\eta$ is a centered Gaussian process defined by
\begin{equation}
\label{E:cov:eta}
  \eta(\mv{w}) = \zeta(\mv{w}) - \sum_{j=1}^p \lambda_j(\mv{w}) \, \zeta(\mv{e}_j), \qquad \mv{w} \in \simplex,
\end{equation}
with $\zeta$ as in Proposition~\ref{P:naive}.
\end{theorem}

\begin{remark}[Covariance function]
\label{R:cov}
The covariance function \eqref{E:cov:zeta} can be expressed in terms of $A$ as follows. An application of the identity $\log(x) = \int_0^\infty \{ \1(s \le x) - \1(s \le 1) \} \, s^{-1} \, ds$ for $x \in (0, \infty)$ yields, by Fubini's theorem,
\begin{multline*}
  \cov \bigl( - \log \xi_i(\mv{v}), \, - \log \xi_i(\mv{w}) \bigr) \\
  \shoveleft{= \int_0^\infty \int_0^\infty 
    \Bigl( P \bigl( \xi_i(\mv{v}) \ge s, \, \xi_i(\mv{w}) \ge t \bigr) - P \bigl( \xi_i(\mv{v}) \ge s \bigr) \, P \bigl( \xi_i(\mv{w}) \ge t \bigr) \Bigr) 
    \frac{ds}{s} \, \frac{dt}{t}} \\
  \shoveleft{= \int_0^\infty \int_0^\infty [ \exp \{ - \ell((w_1s) \vee (v_1t), \ldots, (w_ps) \vee (v_pt)) \}} \\
    - \exp \{ - s \, A(\mv{v}) \} \exp \{ - t \, A(\mv{w}) \} ] \frac{ds}{s} \, \frac{dt}{t}.
\end{multline*}
where $\ell(\mv{y}) = |\mv{y}| \, A(\mv{y}/|\mv{y}|)$ and $|\mv{y}| = |y_1| + \cdots + |y_p|$. Replacing $A$ by any estimator of it results in an estimator of the covariance function. However, a more practical way to estimate this function is by the sample covariance of the pairs $( - \log \xi_i(\mv{v}), \, - \log \xi_i(\mv{w}) )$; see also (the proof of) Theorem~\ref{T:OLS}.
\end{remark}

\begin{remark}[Shape constraints]
A further enhancement to the CFG-estimator is to replace it by the convex minorant of the function
\[
  \min [ \max \{ \Acfg(\mv{w}), w_1, \ldots, w_p \}, 1 ], \qquad \mv{w} \in \simplex,
\]
as in \cite{Deheuvels91} and \cite{JVDF01} for the bivariate case. Although the resulting estimator would be a convex function respecting the bounds $\max(w_1, \ldots, w_p) \le A(\mv{w}) \le 1$, in case $p \ge 3$ this would still not guarantee it to be a genuine Pickands dependence function. Still other ways to impose (some of) the shape restrictions are spline smoothing under constraints \citep{HT00}, orthogonal projection \citep{FGS08}, or Bayesian nonparametrics \citep{GP08}.
\end{remark}

\begin{remark}[Pickands estimator]
A different way to exploit the exponentiality of the random variables $\xi_i(\mv{w})$ in \eqref{E:xiexp} would be via the Pickands estimator
\[
  \frac{1}{\Apick(\mv{w})} = \frac{1}{n} \sum_{i=1}^n \xi_i(\mv{w})
\]
as in \cite{Pickands81}. To impose the vertex constraints $A(\mv{e}_j) = 1$, the techniques of \cite{Deheuvels91} or \cite{HT00} can be used, see \citet[p.~578]{ZWP08}. In the bivariate case however, it is known that the resulting estimators are outperformed by the CFG-estimator $\Acfg$ \citep{S07Ahs, GS09}. This is confirmed in the simulation study in \citet[Section~3]{ZWP08}, as well as by our own simulations in Section~\ref{S:simul}. For this reason, we restrict attention here to the family of CFG-estimators.
\end{remark}

\section{The OLS-estimator}
\label{S:OLS}

The question remains which weight functions $\lambda_j$ to choose in the CFG-estimator \eqref{E:Acfg}. In \cite{ZWP08}, the choice $\lambda_j(\mv{w}) = w_j$ was recommended as a pragmatic one. The option of using variance-minimizing functions $\lambda_j$ was mentioned but not carried out. By casting the estimation problem in a linear regression framework, we will obtain an estimator with the same asymptotic performance as the CFG-estimator with those optimal weights. In this section, we define the estimator and prove its consistency and asymptotic normality, both in the functional sense. In the next section, the gain in efficiency is assessed by means of simulations.

In view of Theorem~\ref{T:CFG}, for each $\mv{w} \in \simplex$ we have
\[
  \sqrt{n} \bigl( \Acfg(\mv{w}) - A(\mv{w}) \bigr) \dto A(\mv{w}) \, \eta(\mv{w}), \qquad n \to \infty,
\]
where $\eta(\mv{w})$ is a zero-mean normal random variable. We will look for those $\lambda_j(\mv{w})$ that minimise the variance of $\eta(\mv{w})$. Let $\zeta$ be the Gaussian process on $\cont(\simplex)$ in Proposition~\ref{P:naive}. For ease of notation, put
\begin{align*}
  \mv{\lambda}(\mv{w})
  &= \bigl( \lambda_1(\mv{w}), \ldots, \lambda_p(\mv{w}) \bigr)^\top, &
  \mv{\zeta}(\mv{e})
  &= \bigl( \zeta(\mv{e}_1), \ldots, \zeta(\mv{e}_p) \bigr)^\top,
\end{align*}
the symbol ``$\top$'' denoting matrix transposition. Then
\begin{align*}
  \var \eta(\mv{w}) 
  &= \var \bigl( \zeta(\mv{w}) - \mv{\lambda}(\mv{w})^\top \, \mv{\zeta}(\mv{e}) \bigr) \\
  &= \var \zeta(\mv{w}) 
  - 2 \, \mv{\lambda}(\mv{w})^\top \, E [\mv{\zeta}(\mv{e}) \, \zeta(\mv{w})]
  + \mv{\lambda}(\mv{w})^\top \, E [ \mv{\zeta}(\mv{e}) \, \mv{\zeta}(\mv{e})^\top ] \, \mv{\lambda}(\mv{w}).
\end{align*}
Note that 
\begin{equation}
\label{E:Sigma}
  \Sigma = E [ \mv{\zeta}(\mv{e}) \, \mv{\zeta}(\mv{e})^\top ]
\end{equation}
is the covariance matrix of $(- \log \xi(\mv{e}_1), \ldots, - \log \xi(\mv{e}_p))^\top$. Provided this matrix is non-singular, $\var \eta(\mv{w})$ attains a unique global minimum for $\mv{\lambda}(\mv{w})$ equal to
\begin{equation}
\label{E:lambdaopt}
  \mv{\lambda}\opt(\mv{w}) = \Sigma^{-1} \, E [\mv{\zeta}(\mv{e}) \, \zeta(\mv{w})].
\end{equation}
With this choice of the weight functions, the variance of 
\begin{equation}
\label{E:etaopt}
  \etaopt(\mv{w}) = \zeta(\mv{w}) - \mv{\lambda}\opt(\mv{w})^\top \, \mv{\zeta}(\mv{e})
\end{equation}
is equal to
\begin{equation}
\label{E:varetaopt}
  \var \etaopt(\mv{w}) 
  = \var \zeta(\mv{w}) - E [\zeta(\mv{w}) \, \mv{\zeta}(\mv{e})^\top] \, \Sigma^{-1} \, E [\mv{\zeta}(\mv{e}) \, \zeta(\mv{w})].
\end{equation}
This variance is minimal over all possible choices of weight functions $\lambda_j$.

The optimal weight functions $\lambda_j\opt$ in \eqref{E:lambdaopt} depend on the unknown Pickands dependence function $A$. Fortunately, replacing these weight functions by uniformly consistent estimators $\hat{\lambda}_{n,j}$ is just as good asymptotically. For such estimated weight functions, define the adaptive CFG-estimator by
\begin{equation}
\label{E:Acfgad}
  \log \Acfgad(\mv{w}) = \log \Anaive(\mv{w}) - \sum_{j=1}^p \hat{\lambda}_{n,j}(\mv{w}) \, \log \Anaive(\mv{e}_j),
\end{equation}

\begin{proposition}[Adaptive CFG-estimator]
\label{P:CFGad}
Assume that, in addition to the assumptions in Proposition~\ref{P:naive}, the matrix $\Sigma$ in \eqref{E:Sigma} is non-singular and $\hat{\lambda}_{n,j}$ are random elements in $\cont(\simplex)$ such that, for every $j \in \{1, \ldots, p\}$,
\[
  \sup_{\mv{w} \in \simplex} | \hat{\lambda}_{n,j}(\mv{w}) - \lambda_j\opt(\mv{w}) | \to 0,
  \qquad n \to \infty, \qquad \text{almost surely},
\]
with $\lambda_j\opt$ as in \eqref{E:lambdaopt}. Then the adaptive CFG-estimator in \eqref{E:Acfgad} satisfies
\begin{equation}
\label{E:CFGad:cons}  
\sup_{\mv{w} \in \simplex} |\Acfgad(\mv{w}) - A(\mv{w})| \to 0, \qquad n \to \infty, \qquad \text{almost surely},
\end{equation}
and in $\cont(\simplex)$,
\begin{equation}
\label{E:CFGad:CLT}
  \sqrt{n} (\Acfgad - A) \dto A \, \etaopt, \qquad n \to \infty,
\end{equation}
where $\etaopt$ is the zero-mean Gaussian process defined in \eqref{E:etaopt}.
\end{proposition}

Finally we propose a particularly convenient way to implement the adaptive CFG-estimator in \eqref{E:Acfgad}. For $\mv{w} \in \simplex$, let $\hat{\mv{\beta}}_n(\mv{w}) = ( \hat{\beta}_{n,0}(\mv{w}), \ldots, \hat{\beta}_{n,p}(\mv{w}) )^\top$ be the minimizer in $(b_0, \ldots, b_p)^\top$ of
\begin{equation}
\label{E:OLS}
  \sum_{i=1}^n \biggl( \bigl(- \log \xi_i(\mv{w}) - \gamma \bigr) - b_0 - \sum_{j=1}^p b_j \bigl(- \log \xi_i(\mv{e}_j) - \gamma \bigr) \biggr)^2.
\end{equation}
In words, $\hat{\mv{\beta}}_n(\mv{w})$ is the ordinary least-squares (OLS) estimator of the vector of regression coefficients in a linear regression of the dependent variable $- \log \xi_i(\mv{w}) - \gamma$ upon the explanatory variables $- \log \xi_i(\mv{e}_j) - \gamma$, $j \in \{1, \ldots, p\}$. Define the OLS-estimator of $A$ via the estimated intercept by
\[
  \log \Aols(\mv{w}) = \hat{\beta}_{n,0}(\mv{w}), \qquad \mv{w} \in \simplex.
\]
Since the residuals
\[
  \hat{\epsilon}_{n,i}(\mv{w}) 
  = \bigl(- \log \xi_i(\mv{w}) - \gamma \bigr) - \hat{\beta}_{n,0}(\mv{w}) - \sum_{j=1}^p \hat{\beta}_{n,j}(\mv{w}) \, \bigl(- \log \xi_i(\mv{e}_j) - \gamma \bigr)
\]
verify $\sum_{i=1}^n \hat{\epsilon}_{n,i}(\mv{w}) = 0$, we have
\begin{equation}
\label{E:Aols}
  \log \Aols(\mv{w}) 
  = \hat{\beta}_{n,0}(\mv{w}) 
  = \log \Anaive(\mv{w}) - \sum_{j=1}^p \hat{\beta}_{n,p}(\mv{w}) \, \log \Anaive(\mv{e}_j),
\end{equation}
that is, the OLS-estimator is equal to the adaptive CFG-estimator with estimated weights $\hat{\lambda}_{n,j}(\mv{w}) = \hat{\beta}_{n,j}(\mv{w})$.  The variance of the (logarithm of the) OLS-estimator can be estimated by the sample variance of the residuals, properly corrected for the loss in number of degrees of freedom,
\begin{equation}
\label{E:sigmaols}
  \sigmaols(\mv{w}) = \frac{1}{n-p-1} \sum_{i=1}^n \hat{\epsilon}_{n,i}^2(\mv{w}), \qquad \mv{w} \in \simplex.
\end{equation}

\begin{theorem}[OLS-estimator]
\label{T:OLS}
Assume that, in addition to the assumptions in Proposition~\ref{P:naive}, the matrix $\Sigma$ in \eqref{E:Sigma} is non-singular. Then, with probability tending to one, the minimizer $\hat{\mv{\beta}}_n(\mv{w})$ of \eqref{E:OLS} is uniquely defined and for $j \in \{1, \ldots, p\}$,
\begin{equation}
\label{E:beta2lambda}
  \sup_{\mv{w} \in \simplex} | \hat{\beta}_{n,j}(\mv{w}) - \lambda_j\opt(\mv{w}) | \to 0,
  \qquad n \to \infty, \qquad \text{almost surely.}
\end{equation}
As a consequence, the OLS-estimator in \eqref{E:Aols} is uniformly consistent,
\begin{equation}
\label{E:Aols:cons}  
\sup_{\mv{w} \in \simplex} |\Aols(\mv{w}) - A(\mv{w})| \to 0, \qquad n \to \infty, \qquad \text{almost surely},
\end{equation}
and in $\cont(\simplex)$,
\begin{equation}
\label{E:Aols:CLT}
  \sqrt{n} (\Aols - A) \dto A \, \etaopt, \qquad n \to \infty,
\end{equation}
where $\etaopt$ is the zero-mean Gaussian process defined in \eqref{E:etaopt}. In addition, the variance estimator in \eqref{E:sigmaols} satisfies
\[
  \sup_{\mv{w} \in \simplex} | \sigmaols(\mv{w}) - \var \etaopt(\mv{w}) | \to 0,
  \qquad n \to \infty, \qquad \text{almost surely}.
\]
\end{theorem}

\begin{remark}[Non-singularity assumption]
In the bivariate case, the assumption that the covariance matrix $\Sigma$ in \eqref{E:Sigma} is non-singular is equivalent to the assumption that the copula $C$ is not the comonotone copula \citep{S07Ahs}. We conjecture that in the general multivariate case, a necessary and sufficient condition for $\Sigma$ to be non-singular is that none of the bivariate margins of $C$ is equal to the comonotone copula.
\end{remark}

\section{Simulations}
\label{S:simul}

In order to investigate the finite-sample properties of the estimators discussed in the previous sections, we generated pseudo-random samples from trivariate extreme-value copulas of logistic type as presented in \cite{tawn1990}:
\begin{multline}
\label{E:logistic}
  A(\mv{w}) = 
  (\theta^{r} w_{1}^{r} + \phi^{r} w_{2}^{r} )^{1/r} + (\theta^{r} w_{2}^{r} + \phi^{r} w_{3}^{r} )^{1/r} + (\theta^{r} w_{3}^{r} + \phi^{r} w_{1}^{r} )^{1/r}  \\
  + \psi ( w_{1}^{r} + w_{2}^{r} + w_{3}^{r}   )^{1/r} + 1 - \theta - \phi -\psi,
  \qquad \mv{w} \in \Delta_p,
\end{multline}
for $(r, \theta, \phi, \psi) \in [1, \infty) \times [0, 1]^3$. To facilitate comparisons, we opted for the same parameter values as chosen in \cite{ZWP08}: a symmetric case, $(r, \theta, \phi, \psi) = (3, 0, 0, 1)$, and an asymmetric one, $(r, \theta, \phi, \psi) = (6, 0.6, 0.3, 0)$. For each case $10\,000$ samples were generated of size $n \in \{50, 100, 200\}$ using the simulation algorithms in \cite{Stephenson03} and implemented in the R-package \textsf{evd} \citep{evd}.

\begin{figure}
\begin{center}
\begin{tabular}{cc}
\includegraphics[width=0.45\textwidth]{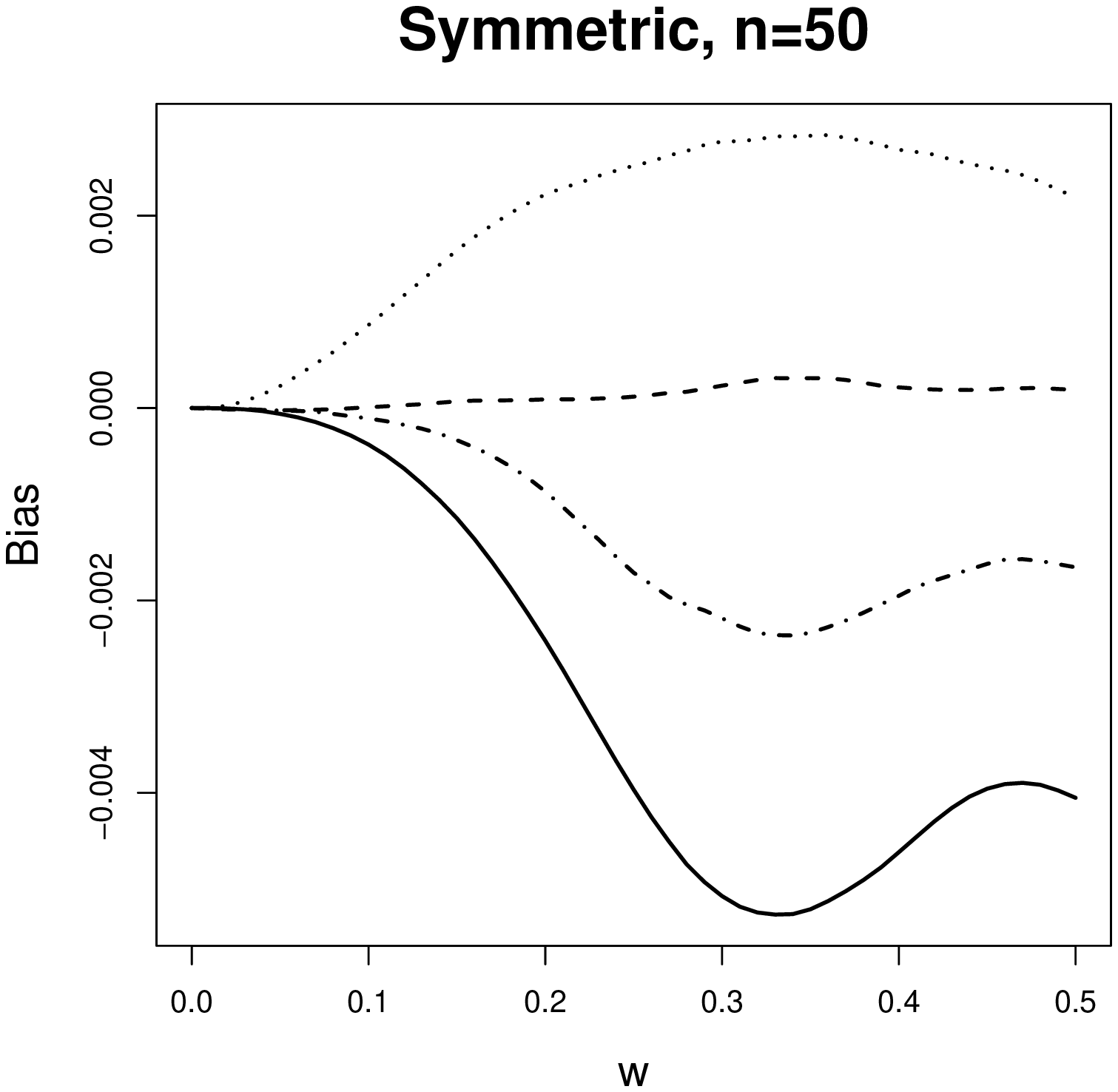}&
\includegraphics[width=0.45\textwidth]{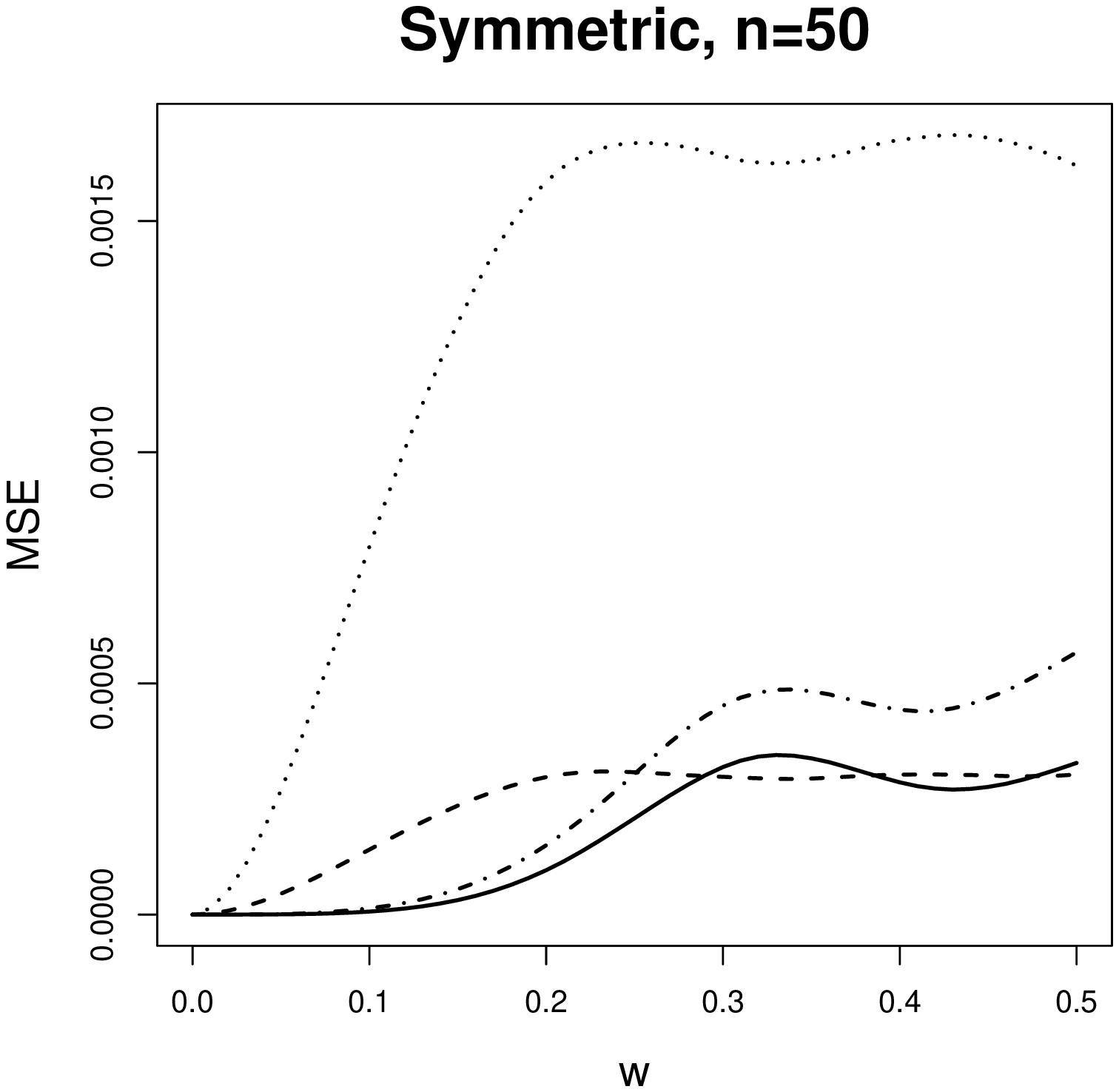}\\
\includegraphics[width=0.45\textwidth]{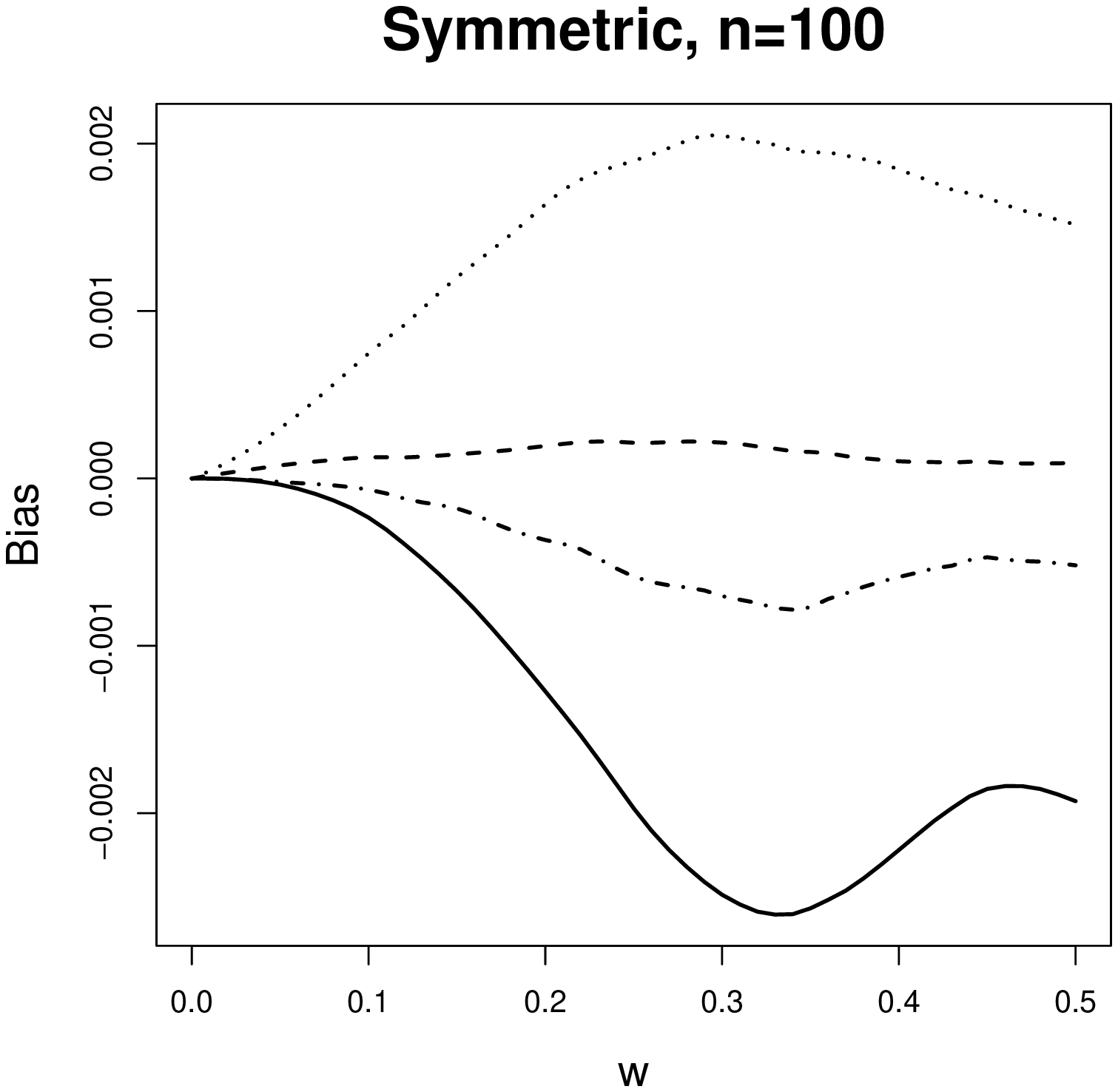}&
\includegraphics[width=0.45\textwidth]{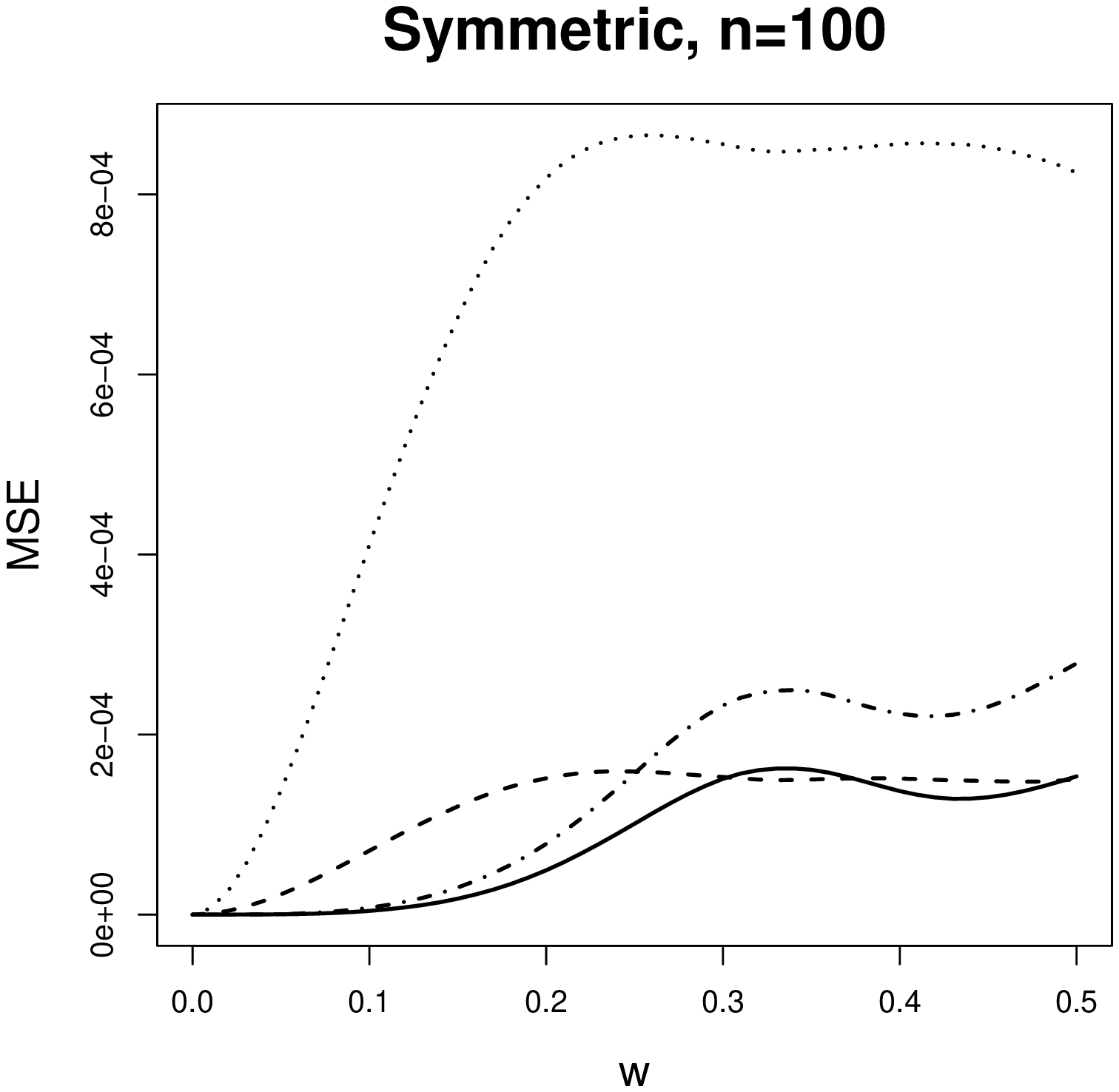}\\
\includegraphics[width=0.45\textwidth]{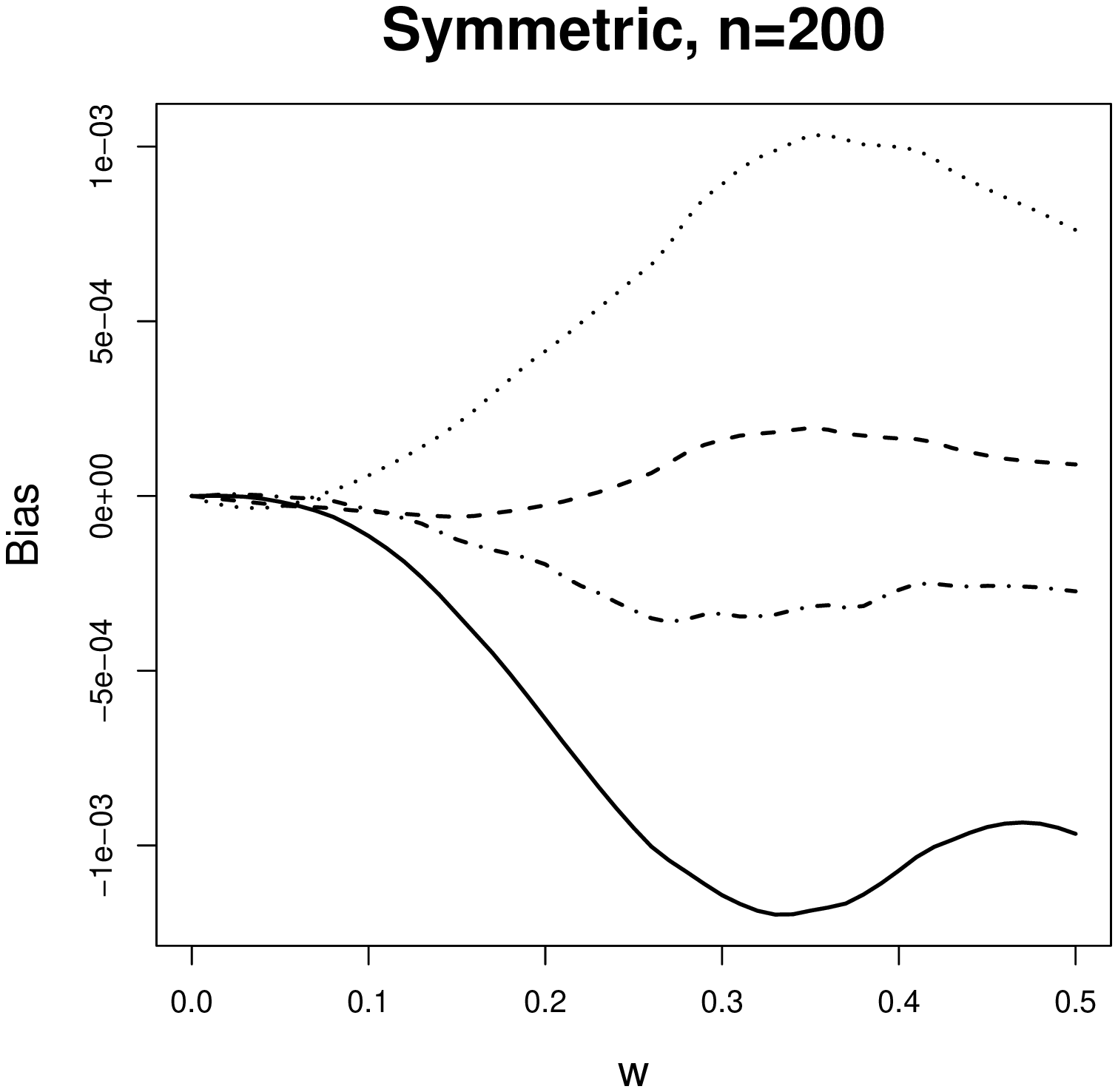}&
\includegraphics[width=0.45\textwidth]{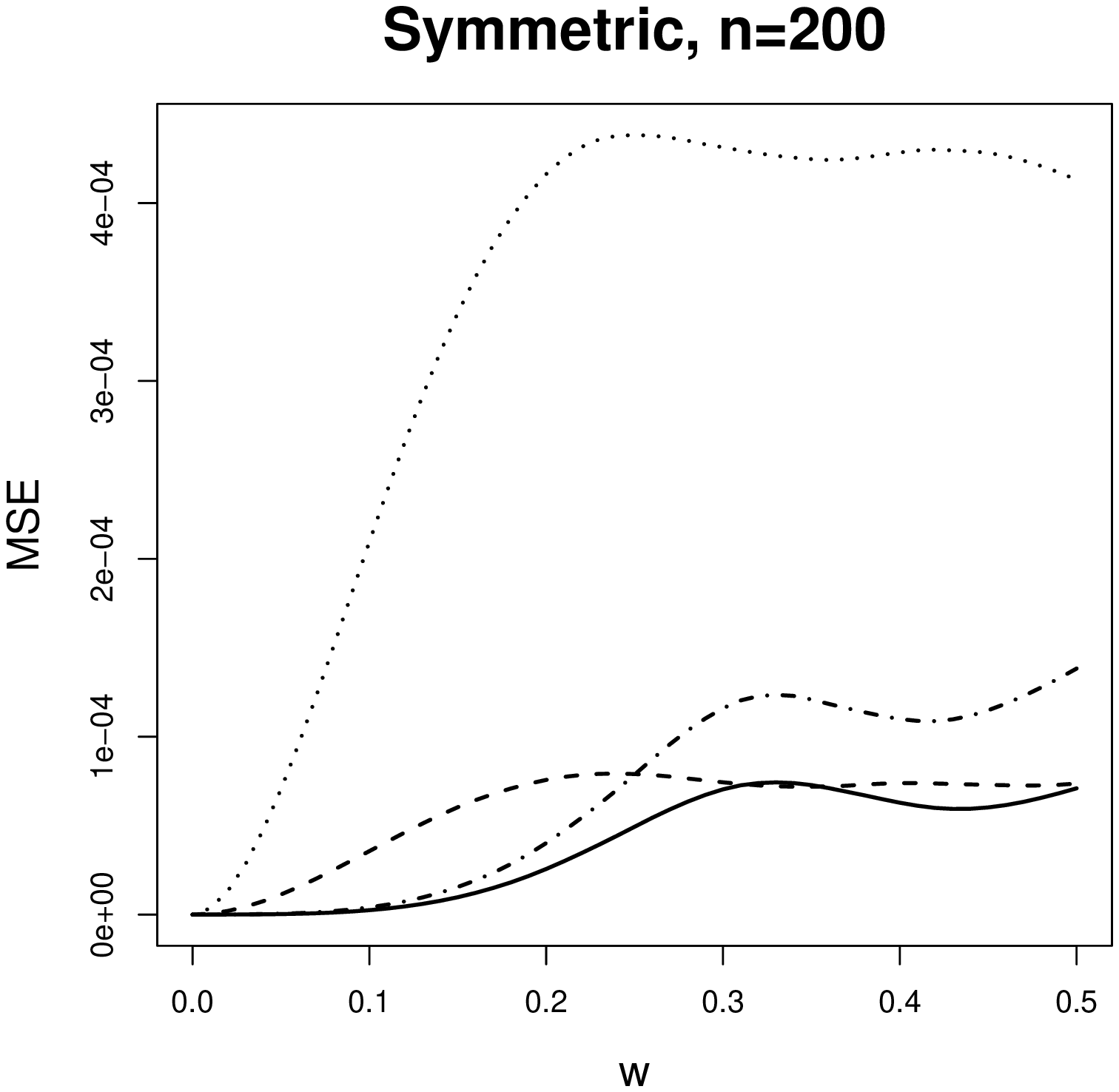}
\end{tabular}
\end{center}
\caption{\label{F:Sym} Biases (left) and mean squared errors (right) of $\Aols(\mv{w})$ (solid), $\Acfg(\mv{w})$ (dashed), $\Aht(\mv{w})$ (dash-dotted) and $\Adh(\mv{w})$ (dotted) along the line $w_1 = w_2$ for $10\,000$ samples of size $n \in \{50, 100, 200\}$ from the trivariate extreme-value copula $C$ with symmetric logistic dependence function $A(\mv{w}) = (w_1^r + w_2^r + w_3^r)^{1/r}$ at $r = 3$.}
\end{figure}

\begin{figure}
\begin{center}
\begin{tabular}{cc}
\includegraphics[width=0.45\textwidth]{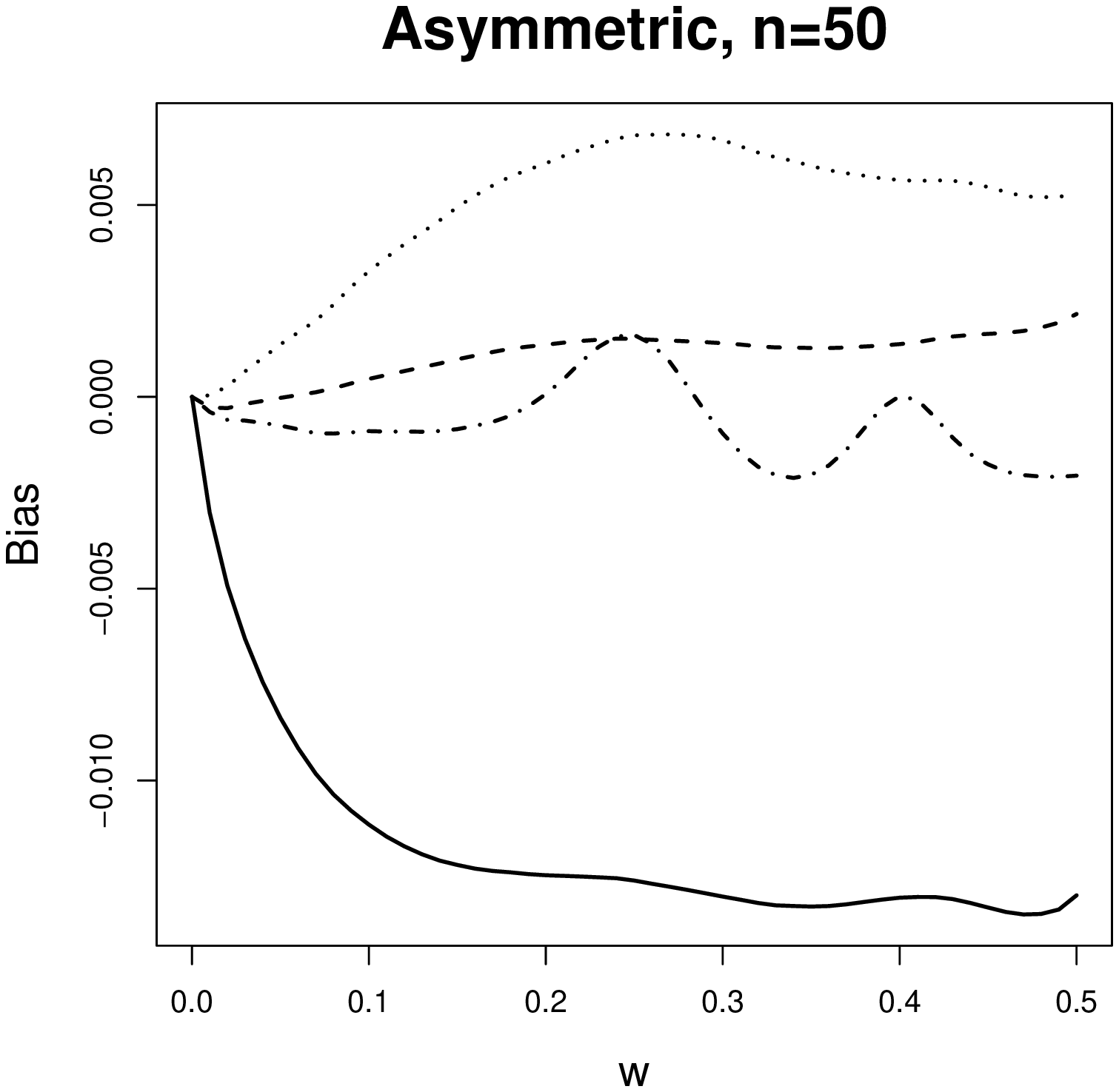}&
\includegraphics[width=0.45\textwidth]{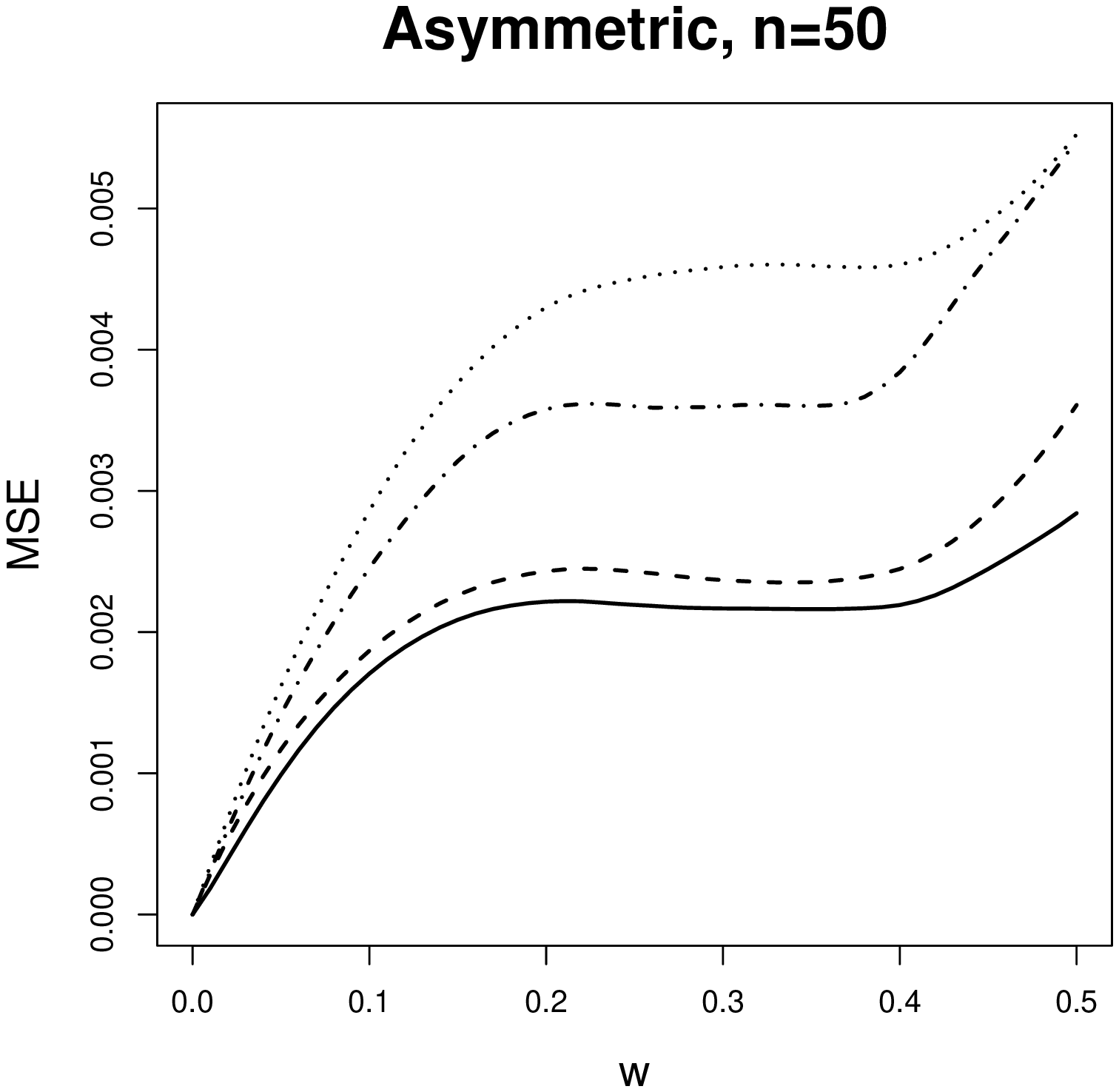}\\
\includegraphics[width=0.45\textwidth]{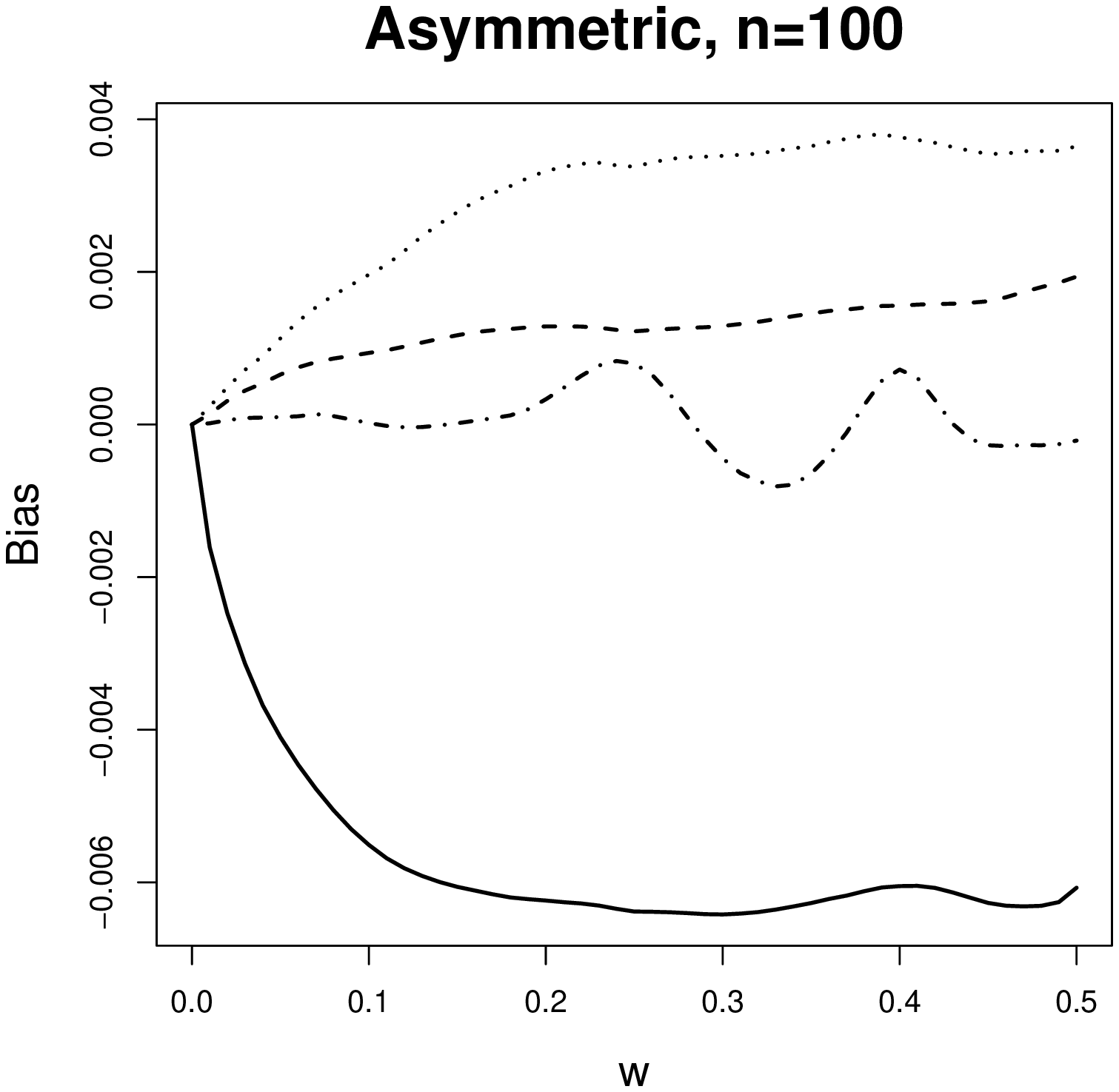}&
\includegraphics[width=0.45\textwidth]{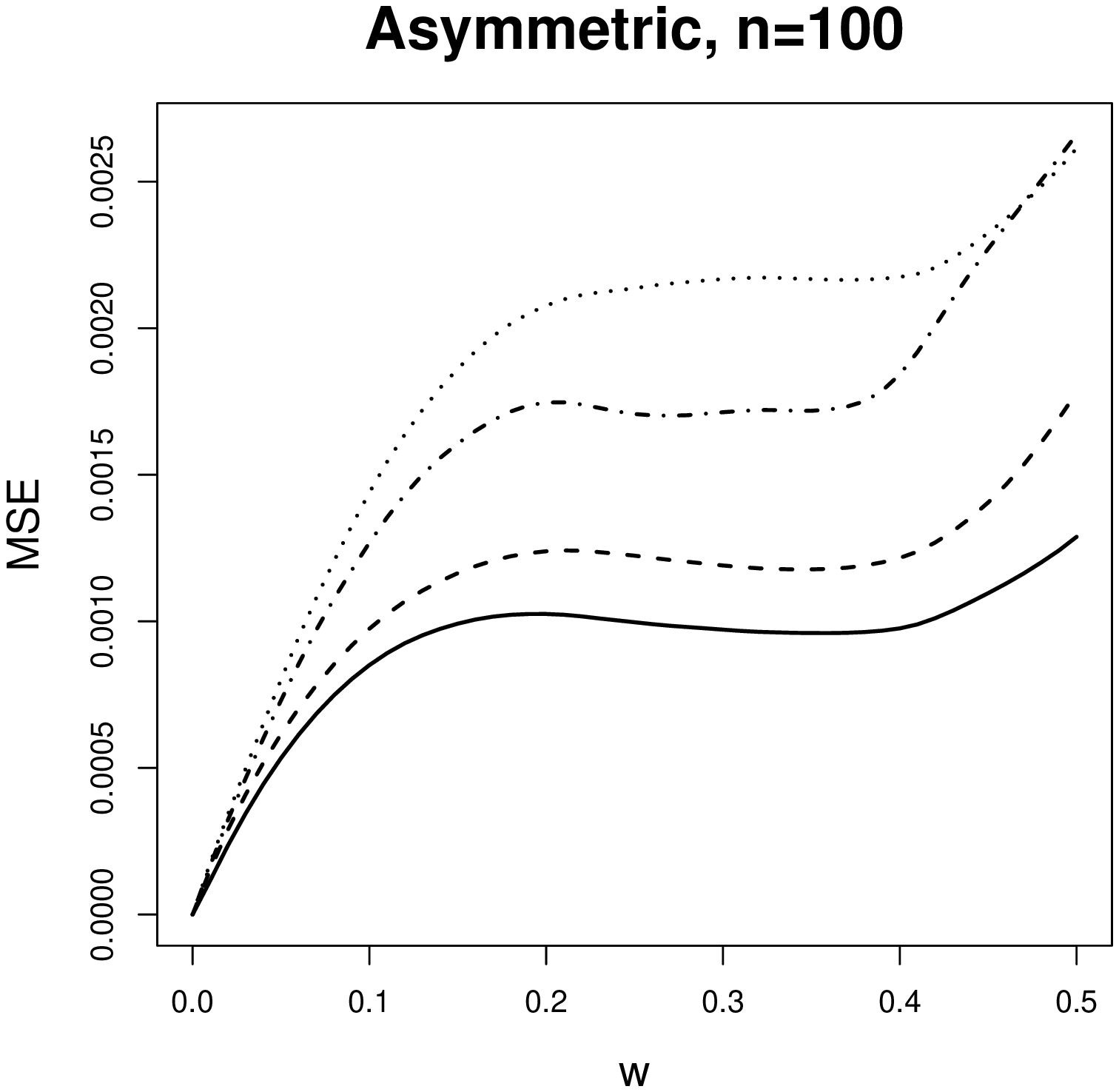}\\
\includegraphics[width=0.45\textwidth]{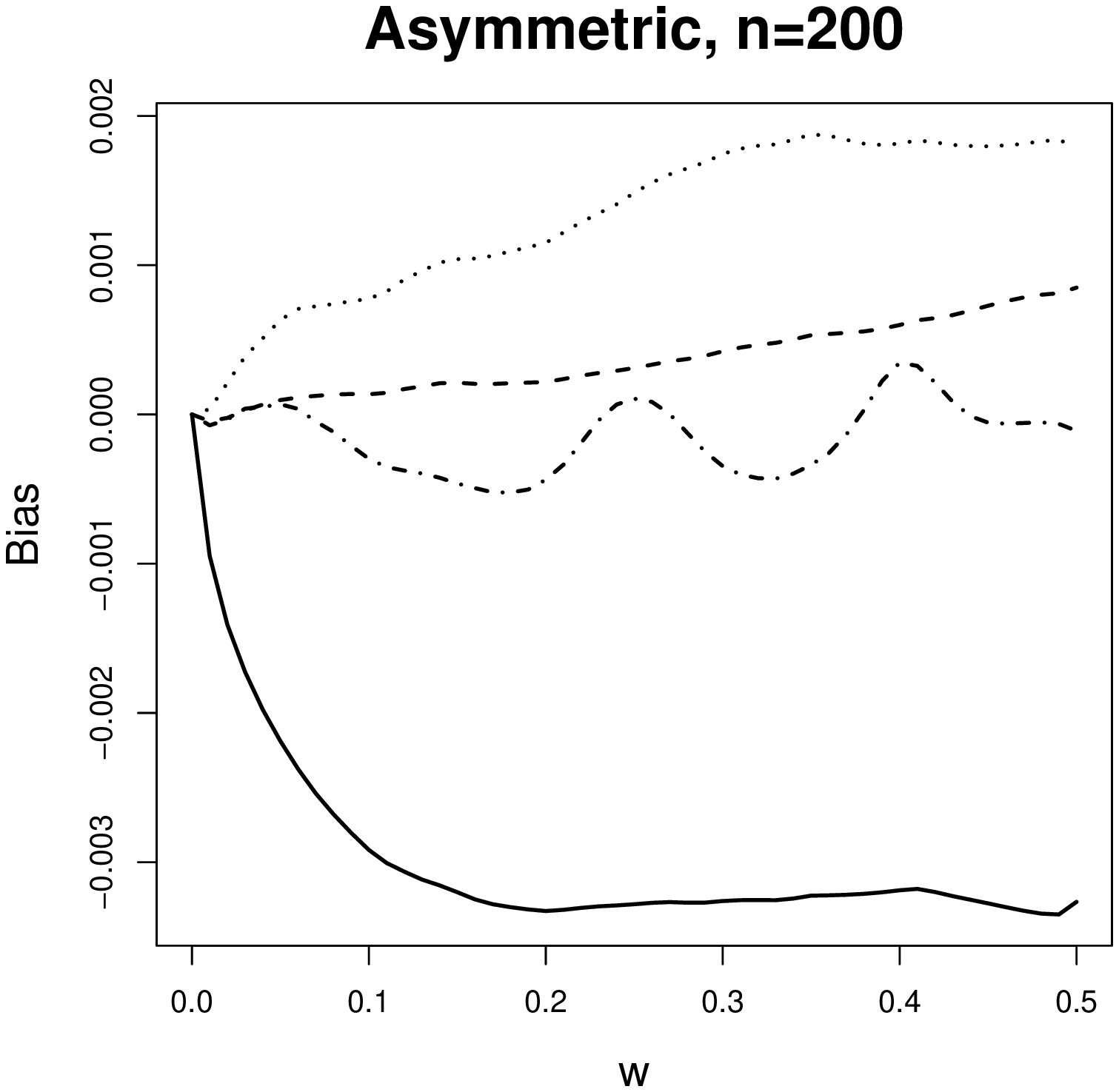}&
\includegraphics[width=0.45\textwidth]{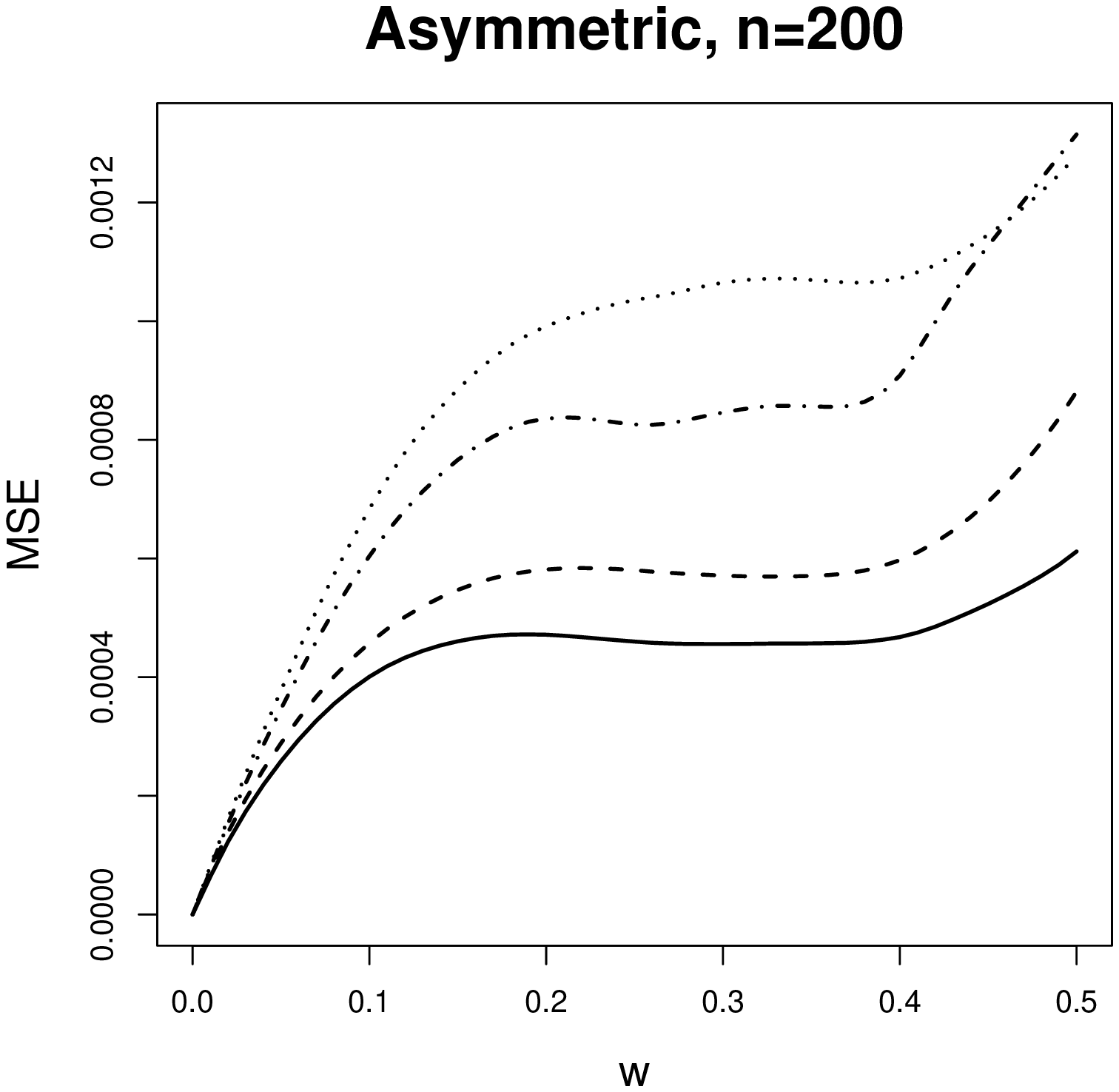}
\end{tabular}
\end{center}
\caption{\label{F:ASym} Biases (left) and mean squared errors (right) of $\Aols(\mv{w})$ (solid), $\Acfg(\mv{w})$ (dashed), $\Aht(\mv{w})$ (dash-dotted) and $\Adh(\mv{w})$ (dotted) along the line $w_1 = w_2$ for $10\,000$ samples of size $n \in \{50, 100, 200\}$ from the trivariate extreme-value copula $C$ with asymmetric logistic dependence function $A$ in \eqref{E:logistic} for $(r, \theta, \phi, \psi) = (6, 0.6, 0.3, 0)$.}
\end{figure}

Four estimators were compared: the CFG-estimator $\Acfg$ with weight functions $\lambda_j(\mv{w}) = w_j$ \citep[as recommended in][]{ZWP08}, the OLS-estimator $\Aols$ in \eqref{E:Aols}, and the enhanced versions of the original Pickands estimator due to \cite{Deheuvels91} and \cite{HT00} as presented in \cite{ZWP08}. To visualize the performances of the estimators, we plotted their biases and mean squared errors along the line $\{ \mv{w} \in \simplex : w_{1} = w_{2} \}$; see Figures~\ref{F:Sym} and~\ref{F:ASym} for the symmetric and asymmetric logistic dependence functions respectively.

In accordance to the theory, the OLS-estimator is in virtually all cases considered more efficient than the CFG-estimator. Moreover, our simulations confirm the findings in \cite{ZWP08} that the CFG-estimator is typically more efficient than the ones of Deheuvels and Hall--Tajvidi. Note that the finite-sample bias of the OLS-estimator is somewhat larger than for the other estimators. However, thanks to its minimum-variance property it ends up as an overall winner in terms of mean squared error.

\appendix
\section{Proofs for Section~\ref{S:CFG}}

\begin{proof}[Proof of Proposition~\ref{P:naive}]
For $\mv{w} \in \simplex$, define $f_{\mv{w}} : (0, \infty)^p \to \RR$ by
\begin{equation}
\label{E:fw}
  f_{\mv{w}}(\mv{y}) = - \log \biggl( \bigwedge_{j=1}^p \frac{y_j}{w_j} \biggr) - \gamma, \qquad \mv{y} \in (0, \infty)^p.
\end{equation}
We can write
\[
  \log \Anaive(\mv{w}) = \frac{1}{n} \sum_{i=1}^n f_{\mv{w}}(\mv{Y}_i).
\]
Consider the function class $\Fclass = \{ f_{\mv{w}} : \mv{w} \in \simplex \}$. We will show that $\Fclass$ is $P$-Donsker and therefore also $P$-Glivenko--Cantelli, where $P$ denotes the common probability distribution on $(0, \infty)^p$ of the random vectors $\mv{Y}_i$. According to Theorem~2.6.8 in \cite{VW96} and the proof thereof, we need to verify that $\Fclass$ is a pointwise separable Vapnik--\u{C}ervonenkis-class (VC-class) that admits an envelope function with a finite second moment under $P$. Pointwise separability follows from the fact that the map $\mv{w} \mapsto f_{\mv{w}}(\mv{y})$ is continuous in $\mv{w} \in \simplex$ for each $\mv{y} \in (0, \infty)^p$. The VC-property can be established by repeated applications of Lemmas~2.6.15 and~2.6.18, items~(i) and~(viii), in \citet{VW96}. Finally, the readily established bound
\begin{equation}
\label{E:fbound}
  \biggl| \log \bigwedge_{j=1}^p \frac{y_j}{w_j} \biggr|
  \le \max \biggl\{ \biggl| \log \bigwedge_{j=1}^p y_j \biggr|, \log (p) + \sum_{j=1}^p | \log y_j | \biggr\}
\end{equation}
yields an envelope function of $\Fclass$ all of whose moments are finite under $P$. Observe that the distribution of $\bigwedge_{j=1}^p Y_{ij}$ is Exponential with mean equal to $\{p \, A(1/p, \ldots, 1/p)\}^{-1} \in [1/p, 1]$.

From the fact that $\Fclass$ is $P$-Glivenko--Cantelli it follows that 
\begin{multline*}
  \sup_{\mv{w} \in \simplex} | \log \Anaive(\mv{w}) - \log A(\mv{w}) | \\
  = \sup_{\mv{w} \in \simplex} \biggl| \frac{1}{n} \sum_{i=1}^n f_{\mv{w}}(\mv{Y}_i) - E [f_{\mv{w}}(\mv{Y})] \biggr|
  \to 0, \qquad n \to \infty, \qquad \text{almost surely.}
\end{multline*}
(Here, we dropped a subscript $i$ for convenience.) Continuity of the map $\exp : \cont(\simplex) \to \cont(\simplex) : f \mapsto \exp(f)$ yields uniform consistency as in \eqref{E:naive:cons}.

Moreover, the $P$-Donsker property entails
\begin{equation}
\label{E:FCLT:ell}
  \sqrt{n} (\log \Acfg - \log A) \dto \zeta, \qquad n \to \infty,
\end{equation}
in the space $\bounded(\simplex)$ of bounded functions from $\simplex$ into $\RR$ equipped with the topology of uniform convergence, where we identified $\Fclass$ with $\simplex$. The process $\zeta$ is zero-mean Gaussian with covariance function given in~\eqref{E:cov:zeta}. The sample paths of the limit process $\zeta$ are continuous with respect to the standard deviation (semi-)metric $\rho$ on $\simplex$ defined by
\[
  \rho(\mv{v}, \mv{w}) = [ \var \{ f_{\mv{v}}(\mv{Y}) - f_{\mv{w}}(\mv{Y}) \} ]^{1/2}, \qquad \mv{v}, \mv{w} \in \simplex.
\]
If $\lim_{n \to \infty} \mv{v}_n = \mv{v}$ in $\simplex$ according to the Euclidean metric, then by continuity of $f_{\mv{w}}(\mv{y})$ in $\mv{w}$ and by uniform integrability, also $\lim_{n \to \infty} \rho(\mv{v}_n, \mv{v}) = 0$. (Uniform integrability is checked by using the bound in \eqref{E:fbound}.) It follows that the trajectories of $\zeta$ are also continuous with respect to the Euclidean metric on $\simplex$, that is, $\zeta$ actually takes its values in $\cont(\simplex)$. As the trajectories of the left-hand side in \eqref{E:FCLT:ell} are continuous too, the convergence in \eqref{E:FCLT:ell} takes place not only $\bounded(\simplex)$ but also in $\cont(\simplex)$.

The convergence in~\eqref{E:naive:CLT} follows from the Hadamard-differentiability of the map $\exp : \cont(\simplex) \to \cont(\simplex) : f \mapsto \exp f$ and the functional delta-method \citep[Section~3.9]{VW96}.
\end{proof}

\begin{proof}[Proof of Theorem~\ref{T:CFG}]
Uniform consistency of $\Acfg$ in \eqref{E:CFG:cons} follows from uniform consistency of $\Anaive$ in \eqref{E:naive:cons} and the fact that the functions $\lambda_j$ are continuous, hence bounded.

To show \eqref{E:CFG:CLT}, define $L : \cont(\simplex) \to \cont(\simplex)$ by
\[
  Lf(\mv{w}) = f(\mv{w}) - \sum_{j=1}^p \lambda_j(\mv{w}) \, f(\mv{e}_j)
\]
for $f \in \cont(\simplex)$ and $\mv{w} \in \simplex$. The operator $L$ is linear and bounded. We have $\log \Acfg = L(\log \Anaive)$. Moreover, as $A(\mv{e}_j) = 1$ for all $j \in \{1, \ldots, p\}$, also $L(\log A) = \log A$. We find
\[
  \sqrt{n} (\log \Acfg - \log A)
  = L \bigl( \sqrt{n} (\log \Anaive - \log A) \bigr)
  \dto L \zeta = \eta, \qquad n \to \infty.
\]
The weak convergence in \eqref{E:CFG:CLT} follows from the functional delta-method \citep[Section~3.9]{VW96}. The representation $\eta = L \zeta$ coincides with \eqref{E:cov:eta}.
\end{proof}

\section{Proofs for Section~\ref{S:OLS}}

\begin{proof}[Proof of Proposition~\ref{P:CFGad}]
If the optimal weight functions $\lambda_j\opt$ were known, we could consider the optimal CFG-estimator
\[
  \log \Acfgopt(\mv{w}) = \log \Anaive(\mv{w}) - \sum_{j=1}^p \lambda_j\opt(\mv{w}) \, \log \Anaive(\mv{e}_j), \qquad \mv{w} \in \simplex.
\]
By Theorem~\ref{T:CFG}, the optimal CFG-estimator is uniformly consistent \eqref{E:CFG:cons} and is asymptotically normal in the sense of \eqref{E:CFG:CLT} with $\eta = \etaopt$. Now
\[
  | \log \Acfgopt(\mv{w}) - \log \Acfgad(\mv{w}) |
  \le \sum_{j=1}^p | \hat{\lambda}_{n,j}(\mv{w}) - \lambda_j\opt(\mv{w}) | \; | \log \Anaive(\mv{e}_j) |.
\]
By uniform consistency of $\hat{\lambda}_{n,j}$ and asymptotic normality of $\sqrt{n} \, \log \Anaive(\mv{e}_j)$, we obtain, as $n \to \infty$,
\begin{align*}
  \sup_{\mv{w} \in \simplex} | \log \Acfgopt(\mv{w}) - \log \Acfgad(\mv{w}) | &\to 0,  \qquad \text{almost surely}, \\
  \sup_{\mv{w} \in \simplex} \sqrt{n} \, | \log \Acfgopt(\mv{w}) - \log \Acfgad(\mv{w}) | &\dto 0.
\end{align*}
As a consequence, the adaptive CFG-estimator is uniformly consistent \eqref{E:CFGad:cons} and asymptotically normal \eqref{E:CFGad:CLT}.
\end{proof}

\begin{proof}[Proof of Theorem~\ref{T:OLS}]
In analogy to the linear regression framework, define the $n \times (p+1)$ matrix
\[
  \mv{X} =
  \begin{pmatrix}
  1 & - \log \xi_{1}(\mv{e}_{1}) - \gamma & \ldots & \log \xi_{1}(\mv{e}_{p}) - \gamma \\
  \ldots & \ldots & \ldots & \ldots \\
  1 & - \log \xi_{n}(\mv{e}_{1}) - \gamma & \ldots & \log \xi_{n}(\mv{e}_{p}) - \gamma
  \end{pmatrix}
\]
and the $n \times 1$ vector
\[
  \mv{Y}(\mv{w}) = \bigl( - \log \xi_1(\mv{w}) - \gamma, \ldots, - \log \xi_n(\mv{w}) - \gamma \bigr)^\top, \qquad \mv{w} \in \simplex.
\]
(No confusion should arise between this $\mv{Y}(\mv{w})$ and the random vectors $\mv{Y}_i$ in \eqref{E:Yi}.) Provided the matrix $\mv{X}^\top \mv{X}$ is non-singular, the OLS-estimator $\hat{\mv{\beta}}_n(\mv{w})$ is given by
\[
  \hat{\mv{\beta}}_n(\mv{w}) = (\mv{X}^\top \mv{X})^{-1} \, \mv{X}^\top \mv{Y}(\mv{w}).
\]
Recall the functions $f_{\mv{w}}$ in \eqref{E:fw}. For $\mv{v}, \mv{w} \in \simplex$, define $g_{\mv{v}, \mv{w}} : (0, \infty)^p \to \RR$ by
\[
  g_{\mv{v}, \mv{w}}(\mv{y}) = f_{\mv{v}}(\mv{y}) \, f_{\mv{w}}(\mv{y}), \qquad \mv{y} \in (0, \infty)^p.
\]
By \eqref{E:fbound} and by Example~2.10.23 in \cite{VW96}, the function class $\{ g_{\mv{v}, \mv{w}} : \mv{v}, \mv{w} \in \simplex \}$ is $P$-Donsker and thus $P$-Glivenko--Cantelli, where $P$ is the common distribution on $(0, \infty)^p$ of the random vectors $\mv{Y}_i$. It follows that, almost surely as $n \to \infty$,
\begin{gather}
\label{E:GC:a}
  \frac{1}{n} \mv{X}^\top \mv{X} \to 
  \begin{pmatrix}
  1 & 0 \\
  0 & \Sigma
  \end{pmatrix}, \\
\label{E:GC:b}
  \sup_{\mv{w} \in \simplex} \biggl| \frac{1}{n} \mv{X}^\top \, \mv{Y}(\mv{w}) -   
  \begin{pmatrix}
  \log A(\mv{w}) \\
  E [\mv{\zeta}(\mv{e}) \zeta(\mv{w})]
  \end{pmatrix} \biggr| \to 0,
\end{gather}
As $\Sigma$ is non-singular, we have
\[
  \begin{pmatrix}
  1 & 0 \\
  0 & \Sigma
  \end{pmatrix}^{-1} =
  \begin{pmatrix}
  1 & 0 \\
  0 & \Sigma^{-1}
  \end{pmatrix},
\]
while $\frac{1}{n} \mv{X}^\top \mv{X}$ is with probability tending to one a non-singular matrix too. We find, almost surely and  uniformly in $\mv{w} \in \simplex$,
\begin{multline*}
  \hat{\mv{\beta}}_n(\mv{w}) = \biggl( \frac{1}{n} \mv{X}^\top \mv{X} \biggr)^{-1} \, \frac{1}{n} \mv{X}^\top \mv{Y}(\mv{w}) \\
  \to 
  \begin{pmatrix}
  1 & 0 \\
  0 & \Sigma^{-1}
  \end{pmatrix}
  \begin{pmatrix}
  \log A(\mv{w}) \\
  E [\mv{\zeta}(\mv{e}) \zeta(\mv{w})]
  \end{pmatrix} 
  = 
  \begin{pmatrix}
  \log A(\mv{w}) \\
  \mv{\lambda}\opt(\mv{w})
  \end{pmatrix},
  \qquad n \to \infty.
\end{multline*}
Equation~\eqref{E:beta2lambda} follows. Proposition~\ref{P:CFGad} and equation \eqref{E:Aols} then yield equations \eqref{E:Aols:cons} and \eqref{E:Aols:CLT}.

Finally, for the estimation of the variance, note that it does not matter asymptotically if we divide by $n$ or by $n-p-1$. Elementary calculations yield
\begin{align*}
  \frac{1}{n} \sum_{i=1}^n \hat{\epsilon}_{n,i}^2(\mv{w}) 
  &= \frac{1}{n} \bigl( \mv{Y}(\mv{w}) - \mv{X} \hat{\mv{\beta}}_n(\mv{w}) \bigr)^\top \, \bigl( \mv{Y}(\mv{w}) - \mv{X} \hat{\mv{\beta}}_n(\mv{w}) \bigr) \\
  &= \frac{1}{n} \mv{Y}(\mv{w})^\top \mv{Y}(\mv{w}) - \biggl( \frac{1}{n} \mv{X}^\top \, \mv{Y}(\mv{w}) \biggr)^\top
  \biggl( \frac{1}{n} \mv{X}^\top \mv{X} \biggr)^{-1} \frac{1}{n} \mv{X}^\top \, \mv{Y}(\mv{w}).
\end{align*}
The Glivenko--Cantelli property yields, almost surely and uniformly in $\mv{w} \in \simplex$,
\begin{multline*}
  \frac{1}{n} \mv{Y}(\mv{w})^\top \mv{Y}(\mv{w})
  = \frac{1}{n} \sum_{i=1}^n \bigl( - \log \xi_i(\mv{w}) - \gamma \bigr)^2 \\
  \to E \bigl[ \bigl( - \log \xi_i(\mv{w}) - \gamma \bigr)^2 \bigr]
  = \var \zeta(\mv{w}) + \bigl( \log A(\mv{w}) \bigr)^2, \qquad n \to \infty.
\end{multline*}
In combination with \eqref{E:GC:a} and \eqref{E:GC:b}, we obtain that $n^{-1} \sum_{i=1}^n \hat{\epsilon}_{n,i}^2(\mv{w})$ converges almost surely and uniformly in $\mv{w} \in \simplex$ to
\begin{multline*}
  \var \zeta(\mv{w}) + \bigl( \log A(\mv{w}) \bigr)^2
  - \begin{pmatrix}
  \log A(\mv{w}) \\
  E [\mv{\zeta}(\mv{e}) \zeta(\mv{w})]
  \end{pmatrix}^\top 
    \begin{pmatrix}
  1 & 0 \\
  0 & \Sigma^{-1}
  \end{pmatrix}
  \begin{pmatrix}
  \log A(\mv{w}) \\
  E [\mv{\zeta}(\mv{e}) \zeta(\mv{w})]
  \end{pmatrix} \\
  = \var \zeta(\mv{w}) - E [\mv{\zeta}(\mv{e})^\top \zeta(\mv{w})] \, \Sigma^{-1} \, E [\mv{\zeta}(\mv{e}) \zeta(\mv{w})],
\end{multline*}
which by \eqref{E:varetaopt} is equal to $\var \etaopt(\mv{w})$.
\end{proof}


\end{document}